\documentclass[11pt,a4paper]{article}

\usepackage[T1]{fontenc}
\usepackage[utf8]{inputenc}
\usepackage{lmodern}
\usepackage[margin=1in]{geometry}
\usepackage{amsmath,amssymb,amsfonts,amsthm}
\usepackage{microtype}
\usepackage{authblk}
\usepackage{hyperref}
\hypersetup{
  hidelinks,
  pdftitle={Projective flatness of cylindrically symmetric Finsler metrics},
  pdfauthor={Newton Solorzano, Victor Leon, Dik Lujerio},
  pdfsubject={Differential Geometry; MSC 53B40, 53C60},
  pdfkeywords={Finsler metric, cylindrically symmetric, weakly orthogonally invariant, projectively flat, strong convexity, Hamel equations}
}
\theoremstyle{plain}
\newtheorem{theorem}{Theorem}
\newtheorem{proposition}[theorem]{Proposition}
\newtheorem{lemma}[theorem]{Lemma}
\newtheorem{corollary}[theorem]{Corollary}
\theoremstyle{definition}
\newtheorem{example}{Example}
\theoremstyle{remark}
\newtheorem{remark}{Remark}

\newcommand{\R}{\mathbb R}
\newcommand{\OO}{\mathrm O}
\newcommand{\dd}{\,\mathrm d}
\newcommand{\ip}[2]{\langle #1,#2\rangle}
\newcommand{\Apot}{\mathcal A}
\DeclareMathOperator{\erf}{erf}
\title{Projective flatness of cylindrically symmetric Finsler metrics}
\author[1]{Newton Sol\'orzano}
\author[1]{V\'ictor Le\'on}
\author[2]{Dik Lujerio}
\affil[1]{ILACVN--CICN, Universidade Federal da Integra\c{c}\~ao
Latino-Americana--UNILA, Itaipu Parquetec, Foz do Igua\c{c}u--PR,
85867-970, Brazil}
\affil[2]{Departamento Acad\'emico de Matem\'atica de la Facultad
 de Ciencias, Universidad Nacional Santiago Ant\'unez de Mayolo--UNASAM,
 Av. Universitaria s/n, Shancay\'an, distrito de Independencia,
 Huaraz, Peru}
\affil[ ]{\href{mailto:nmayer159@gmail.com}{\texttt{nmayer159@gmail.com}}\quad
\href{mailto:vicaml19@gmail.com}{\texttt{vicaml19@gmail.com}}\quad
\href{mailto:dlujeriog@unasam.edu.pe}{\texttt{dlujeriog@unasam.edu.pe}}}

\date{}

\begin{document}
\maketitle

\begin{abstract}
In this paper, we study cylindrically symmetric projectively flat Finsler metrics \(F=|\bar y|\phi(x^0,r,s,z)\). We reformulate the conditions for projective flatness and strong convexity in terms of \(\Omega=\phi-s\phi_s-z\phi_z\), isolating the compatibility condition missing in earlier work. Using the method of characteristics, we give an integral formula for all smooth local solutions near points with \(r>0\) and \(s=z=0\). These solutions are determined by a characteristic potential up to the addition of an invariant exact one-form. In the reversible case, the equation for the potential reduces to a \((2+1)\)-dimensional wave equation on a half-space. We illustrate the reconstruction with a polynomial potential and compute the potentials of a known family of spherically symmetric metrics and their linear pullbacks. The convexity conditions also admit decreasing profiles and recover the spherical Beltrami metric as a particular case. Finally, suitable positive sums of linear pullbacks yield smooth reversible projectively flat Finsler metrics on \(\mathbb R^{n+1}\) whose tangent norm at the origin is non-Riemannian. Consequently, these metrics are not linear pullbacks of spherically symmetric metrics.

\end{abstract}



\section{Introduction}\label{sec:introduction}

A Finsler metric on an open subset of $\R^{n+1}$ is called
projectively flat in the given coordinates if its unparameterized
geodesics are straight line segments. Hamel's equations characterize
such metrics:
\begin{equation}\label{eq:Hamel}
 F_{x^Ay^B}y^A-F_{x^B}=0.
\end{equation}
For the definitions and basic properties of Finsler metrics,
we refer to \cite{BCS}.

Let $M\subset\R\times\R^n$, $n\geq2$, be an open domain invariant
under the action $(x^0,\bar x)\mapsto(x^0,O\bar x)$,
$O\in\OO(n)$. A Finsler metric $F$ on $M$ is called
cylindrically symmetric if
\begin{equation}\label{eq:1}
 F\bigl((x^0,O\bar x),(y^0,O\bar y)\bigr)
 =F\bigl((x^0,\bar x),(y^0,\bar y)\bigr)
\end{equation}
for every $O\in\OO(n)$.
These metrics are also called weakly orthogonally invariant
\cite{Liu2024}. This class contains spherically symmetric metrics
and several classes of warped product metrics.

On the region $|\bar y|>0$, every such metric can be written as
\cite{Solorzano2022}
\begin{equation}\label{def:F}
 F(x,y)=|\bar y|\phi(x^0,r,s,z),\qquad
 r=|\bar x|,\quad s=\frac{\ip{\bar x}{\bar y}}{|\bar y|},
 \quad z=\frac{y^0}{|\bar y|}.
\end{equation}
Here, $|\cdot|$ and $\ip{\cdot}{\cdot}$ denote the standard
Euclidean norm and inner product on $\R^n$.
Sol\'orzano and Le\'on \cite{Solorzano2023} derived a system of
partial differential equations for projective flatness in this class
and obtained an explicit separated family of solutions.
In \cite{Solorzano2023} it was stated that projective flatness is
equivalent to the first two equations of \eqref{eq:PF-omega} below.
Liu, Mo, and Zhu \cite{Liu2024} showed that these two equations are only necessary: adding a non-closed invariant one-form to a projectively flat metric preserves them but destroys projective flatness \cite[Section~3]{Liu2024}. They proved that projective
flatness is equivalent to the system \eqref{eq:corrected}--\eqref{eq:p}
below \cite[Theorem~1.1]{Liu2024}; their function $p_1$ is $D\phi$ in
our notation. Proposition~\ref{prop:equivalent-systems} makes the
difference explicit: the missing condition is
$\phi_{rz}=r\phi_{x^0s}$, which is precisely the closedness condition
for invariant one-forms $a(x^0,r)y^0+b(x^0,r)\ip{\bar x}{\bar y}$
\cite[Lemma~3.3]{Liu2024}.
In \cite[Section~4]{Liu2024}, the same authors obtained projectively
flat weakly orthogonally invariant metrics on $\R^{n+1}$ by pulling
back Bryant's two-parameter metrics on $S^{n+1}$ and adding closed
invariant one-forms, and computed their flag curvature; these metrics
are not spherically symmetric when the two parameters differ (see the
proof of \cite[Proposition~2.2]{Liu2024}).
Our approach is complementary: instead of starting from specific
metrics, we describe all local solutions through a single
characteristic potential.
The spray coefficients and the corresponding criterion were also
considered in \cite{Weakly2026} in the study of metrics with
vanishing Douglas curvature.

Throughout the paper we use the first-order operator
\begin{equation}\label{eq:D}
 D=z\partial_{x^0}+\frac{s}{r}\partial_r+\partial_s .
\end{equation}
In the notation used here, the projective-flatness equations are
\begin{equation}\label{eq:corrected}
 p_2=p_3=0,
\end{equation}
where
\begin{equation}\label{eq:p}
 p_2=(D\phi)_z-2\phi_{x^0},\qquad
 p_3=(D\phi)_s-\frac2r\phi_r.
\end{equation}
In this paper, we study this system through the function
\[
 \Omega=\phi-s\phi_s-z\phi_z.
\]

Our main results are the following.
First, \eqref{eq:corrected} is equivalent to
\[
 \Omega_{x^0}=\phi_{sz},\qquad
 \Omega_r=r\phi_{ss},\qquad
 \phi_{rz}=r\phi_{x^0s}
\]
(Proposition~\ref{prop:equivalent-systems}).
These equations imply $D\Omega=0$, that is,
\[
 z\Omega_{x^0}+\frac{s}{r}\Omega_r+\Omega_s=0,
\]
so $D$ is a characteristic vector field for $\Omega$.
(This transport equation also follows by taking $U=L=0$ in
\cite[Corollary~3]{Weakly2026}.)
Integrating along the characteristic curves, we obtain
\[
 \Omega=\Apot(\xi,q,z),\qquad
 \xi=x^0-zs,\qquad q=r^2-s^2,
\]
where the \emph{characteristic potential} $\Apot$ satisfies
\[
 2\Apot_{qz}+z\Apot_{\xi\xi}=0
\]
(Lemma~\ref{lem:PF-A}). For reversible metrics, $\Apot$ is even in
$z$, and in the variables $(\xi,q,z^2)$ this equation becomes a
$(2+1)$-dimensional wave equation (Remark~\ref{rem:wave}).

Conversely, near a point $(x^0_*,r_*,0,0)$ with $r_*>0$, we give an
integral formula $\Phi_{\Apot}$ for every smooth potential
satisfying the equation above and $\Apot_z(\xi,q,0)=0$, and we show
that every smooth local solution has the form
\[
 \phi=\Phi_{\Apot}+z\psi_{x^0}(x^0,r)+\frac{s}{r}\psi_r(x^0,r),
\]
where $\psi$ is a smooth function on a sufficiently small,
contractible base neighborhood
(Theorem~\ref{thm:PF-reconstruction}).
The additional term corresponds to an invariant exact one-form
evaluated on the tangent vector.
We also express the convexity quantity $\Lambda$ of
Lemma~\ref{lem:convexity} in terms of $\Apot$, which gives a
necessary and sufficient condition for strong convexity on the
regular reduced region (Proposition~\ref{prop:potential-positive}).
Example~\ref{ex:local} applies the reconstruction to a polynomial
potential that does not come from the spherically symmetric family.

The reconstruction solves Hamel's equations locally on the regular
reduced region; positivity, strong convexity, and extension to the
whole slit tangent bundle must be checked separately.
In Section~\ref{sec:construction} we do this for the
projectively flat spherically symmetric metrics classified in
\cite{Zhou2012,HuangMo2013} and for their anisotropic linear
pullbacks $(x^0,\bar x)\mapsto(cx^0,\bar x)$.
For $c=1$ these metrics are known, and for $c\neq1$ they are
linearly equivalent to them. What is new is the computation of their
characteristic potentials and an exact convexity criterion,
$f>0$ and $f+2wf'>0$, which is weaker than monotonicity of the
profile $f$ and also covers the round sphere
(Theorem~\ref{thm:construction}).
Genuinely cylindrical global examples are obtained in
Section~\ref{sec:examples} as positive sums of pullbacks with
different anisotropies (Example~\ref{ex:sum}); they are not linear
pullbacks of spherically symmetric metrics. This gives an elementary
source of non-spherical examples, different from the Bryant-type
metrics of \cite{Liu2024}, together with their characteristic
potentials.

In Section~\ref{sec:preliminaries}, we give the notation and
formulas used in the proofs.
In Section~\ref{sec:PF}, we prove the potential representation
and the reconstruction theorem.
Section~\ref{sec:construction} treats the spherical family and its
pullbacks, and Section~\ref{sec:examples} gives explicit examples.

\section{Preliminaries}\label{sec:preliminaries}

In this section, we give the notation and formulas used in the
proofs. Let $M$ be an $\OO(n)$-invariant open domain in
$\R\times\R^n$, $n\geq2$, with the action given in
Section~\ref{sec:introduction}.
Write
\begin{align}
 x&=(x^0,\bar x),\qquad \bar x=(x^1,\ldots,x^n),\label{coordx}\\
 y&=(y^0,\bar y),\qquad \bar y=(y^1,\ldots,y^n).\label{coordy}
\end{align}
Throughout the paper, the convention for indices is
\[
 0\leq A,B,C\leq n,\qquad 1\leq i,j\leq n.
\]
Repeated indices are summed. We set
\[
 0_M=\{(x,0):x\in M\},\qquad TM_0=TM\setminus0_M.
\]
A Finsler metric is smooth and positive on $TM_0$, positively
homogeneous of degree one in $y$, with positive definite
fundamental tensor
\[
 g_{AB}=\frac12(F^2)_{y^Ay^B}.
\]
For \eqref{def:F}, define
\begin{align}
 \Omega&=\phi-s\phi_s-z\phi_z,\label{DefOmega}\\
 \Lambda&=\Omega\phi_{zz}
 +(r^2-s^2)(\phi_{ss}\phi_{zz}-\phi_{sz}^2).
 \label{eq:Lambda}
\end{align}
The variables $r,s,z$ are used on $|\bar y|>0$. The regular
reduced region is $r>0$, $|s|<r$. Full smoothness at the remaining
directions must be verified in the original variables.
We write $u=|\bar y|$ and $v=\bar y/u$, with components $v^i=y^i/u$.

\begin{lemma}\label{lem:convexity}
Suppose that $\phi>0$. On the regular reduced region,
$F=|\bar y|\phi$ is strongly convex if and only if
\[
 \phi_{zz}>0\qquad\text{and}\qquad \Lambda>0,
\]
together with $\Omega>0$ when $n\geq3$.
\end{lemma}
Lemma~\ref{lem:convexity} is essentially
\cite[Theorem~1]{Solorzano2023}, as used in
\cite[Lemma~4.4]{Liu2024}. We include a short proof, which also makes
the role of the condition $\phi_{zz}>0$ explicit.
\begin{proof}
For fixed $x$, let $d_yF$ and $D_y^2F$ denote the differential
and Hessian of $F(x,\cdot)$ at $y$, respectively. In coordinates,
\[
 d_yF(h)=F_{y^A}h^A,\qquad
 D_y^2F(h,k)=F_{y^Ay^B}h^Ak^B.
\]

For a fiber variation $\delta y$, define
\[
 \omega=\delta\bar y-\ip{v}{\delta\bar y}v,\qquad
 \mu=\delta y^0-z\ip{v}{\delta\bar y},\qquad
 \nu=\ip{\bar x-sv}{\delta\bar y}.
\]
Direct differentiation gives
\[
 D_y^2F(\delta y,\delta y)=\frac1u
 \left(
 \phi_{zz}\mu^2+2\phi_{sz}\mu\nu
 +\phi_{ss}\nu^2+\Omega|\omega|^2
 \right).
\]
Write $q=r^2-s^2>0$. Since $\bar x-sv\perp v$ and
$|\bar x-sv|^2=q$, we may decompose
\[
 \omega=\frac{\nu}{q}(\bar x-sv)+\omega_\perp,
\]
where $\omega_\perp$ is orthogonal to both $v$ and
$\bar x-sv$. Then
\[
 uD_y^2F(\delta y,\delta y)=
 \bigl(\mu,\ \nu/\sqrt q\bigr)
 \begin{pmatrix}
 \phi_{zz}&\sqrt q\,\phi_{sz}\\
 \sqrt q\,\phi_{sz}&\Omega+q\phi_{ss}
 \end{pmatrix}
 \begin{pmatrix}\mu\\ \nu/\sqrt q\end{pmatrix}
 +\Omega|\omega_\perp|^2.
\]
The determinant of the matrix is $\Lambda$.
The linear map
\[
 \delta y\longmapsto(\mu,\nu/\sqrt q,\omega_\perp)
\]
is onto $\R^2\times\R^{n-2}$, identifying the orthogonal
complement with $\R^{n-2}$, and its kernel is spanned by $y$.
By Sylvester's criterion, $D_y^2F$ is positive semidefinite
with kernel $\R y$ if and only if
\[
 \phi_{zz}>0,\qquad \Lambda>0,
\]
together with $\Omega>0$ when $n\geq3$.

Finally, differentiating $F^2$ gives
\[
 g_y=d_yF\otimes d_yF+F D_y^2F.
\]
Homogeneity yields
\[
 d_yF(y)=F,\qquad D_y^2F(y,\cdot)=0.
\]
Since $F>0$, the tensor $g_y$ is positive definite if and only
if $D_y^2F$ is positive semidefinite with kernel $\R y$.
Indeed, if $g_y$ is positive definite, then
$D_y^2F=g_y/F$ on $\ker d_yF$, which is complementary to
$\R y$. Together with the vanishing on $\R y$, this gives
the asserted property of $D_y^2F$.

Conversely, if $D_y^2F\geq0$ with kernel $\R y$, then
$g_y\geq0$. Moreover, $g_y(V,V)=0$ forces $V\in\R y$
and $d_yF(V)=0$. Since $d_yF(y)=F>0$, we obtain $V=0$.
Thus $g_y$ is positive definite.
\end{proof}

The spray coefficients satisfy
\begin{equation}\label{eq:spray-decomp}
 G^A=Py^A+Q^A,\qquad
 P=\frac{F_{x^C}y^C}{2F},\qquad
 Q^A=\frac F2g^{AB}(F_{x^Cy^B}y^C-F_{x^B}).
\end{equation}
We recall the reduced spray formulas from \cite{Weakly2026},
using the notation in \eqref{eq:p}. A proof is included for
the computations used below.

\begin{proposition}\label{prop:spray}
Let $F=|\bar y|\phi(x^0,r,s,z)$ be a smooth Finsler metric.
On the regular reduced region, its spray coefficients are
\[
 G^0=u^2N,\qquad G^i=u^2(Wv^i+Ux^i),
\]
where
\begin{align}
 N&=z(W+sU)+L,\label{def:N}\\
 W&=\frac1\phi\left\{\frac{D\phi}{2}
       -s\phi U-\phi_zL-(r^2-s^2)\phi_sU\right\},\label{def:W}\\
 L&=\frac{(\Omega+(r^2-s^2)\phi_{ss})p_2
                  -(r^2-s^2)\phi_{sz}p_3}{2\Lambda},\label{def:L}\\
 U&=\frac{\phi_{zz}p_3-\phi_{sz}p_2}{2\Lambda}.\label{def:U}
\end{align}
In particular, $F$ is projectively flat if and only if $L=U=0$.
In that case,
\[
 G^A=Py^A,\qquad P=\frac{u\,D\phi}{2\phi}.
\]
\end{proposition}
\begin{proof}
Define the Hamel residuals by
\[
 H_B=F_{x^Ay^B}y^A-F_{x^B}.
\]
Using $F_{x^A}y^A=u^2D\phi$, direct differentiation gives
\[
 H_0=up_2,\qquad
 H_i=u\left[p_3x^i-(zp_2+sp_3)v^i\right].
\]
By invariance and homogeneity, write
\[
 (G^0,\bar G)
 =u^2\bigl(z(W+sU)+L,\,
 (W+sU)v+U(\bar x-sv)\bigr).
\]
The Euler--Lagrange equations and homogeneity give
\[
 2F_{y^Ay^B}G^A=H_B,\qquad
 F_{y^A}G^A=\frac12F_{x^A}y^A.
\]
Using the fiber Hessian formula from the preceding lemma,
the nonradial components of the first identity yield
\[
 \begin{pmatrix}
 \phi_{zz}&q\phi_{sz}\\
 \phi_{sz}&\Omega+q\phi_{ss}
 \end{pmatrix}
 \begin{pmatrix}L\\U\end{pmatrix}
 =\frac12\begin{pmatrix}p_2\\p_3\end{pmatrix},
 \qquad q=r^2-s^2.
\]
Its determinant is $\Lambda>0$, so solving the system
gives \eqref{def:L} and \eqref{def:U}.
The second identity gives
\[
 \phi W+(s\phi+q\phi_s)U+\phi_zL=\frac{D\phi}{2},
\]
which proves \eqref{def:W}; \eqref{def:N} follows from
the spray decomposition.

The invertibility of the system shows that $L=U=0$
if and only if $p_2=p_3=0$. By Hamel's equations,
this is equivalent to projective flatness. In that case,
\[
 G^A=Py^A,\qquad P=uW=\frac{uD\phi}{2\phi}.
\]
Finally, smoothness and density of the regular reduced
region extend Hamel's equations to the remaining
tangent directions by continuity.
\end{proof}

\section{Projectively flat metrics}\label{sec:PF}

\begin{proposition}\label{prop:equivalent-systems}
For a smooth reduced function $\phi$, the equations $p_2=p_3=0$
are equivalent to
\begin{equation}\label{eq:PF-omega}
 \Omega_{x^0}=\phi_{sz},\qquad
 \Omega_r=r\phi_{ss},\qquad
 \phi_{rz}=r\phi_{x^0s}.
\end{equation}
\end{proposition}
\begin{proof}
If $p_2=p_3=0$, equality of the mixed derivatives of $D\phi$ gives
\[
 0=(p_2)_s-(p_3)_z=-2\phi_{x^0s}+\frac2r\phi_{rz}.
\]
Expanding \eqref{eq:p} and substituting this identity,
we obtain
\[
 0=p_2=\phi_{sz}-\Omega_{x^0},\qquad
 0=rp_3=r\phi_{ss}-\Omega_r.
\]
Conversely, expansion gives
\begin{align*}
 p_2&=\phi_{sz}-\Omega_{x^0}
       +s\left(\frac{\phi_{rz}}r-\phi_{x^0s}\right),\\
 p_3&=\phi_{ss}-\frac{\Omega_r}r
       -z\left(\frac{\phi_{rz}}r-\phi_{x^0s}\right).
\end{align*}
Thus \eqref{eq:PF-omega} implies $p_2=p_3=0$.
\end{proof}

\begin{remark}\label{rem:missing}
The first two equations in \eqref{eq:PF-omega} are the system
\cite[(3.13)]{Liu2024}, which was claimed in \cite{Solorzano2023} to
characterize projective flatness. The third equation,
$\phi_{rz}=r\phi_{x^0s}$, is the extra condition. It cannot be
dropped: if $\phi$ is a solution and $a,b$ are smooth functions of
$(x^0,r)$, then $\phi+az+bs$ has the same $\Omega$, $\phi_{ss}$ and
$\phi_{sz}$, so it satisfies the first two equations, while the third
holds if and only if $a_r=rb_{x^0}$, that is, if and only if the
one-form $ay^0+b\ip{\bar x}{\bar y}$ is closed
\cite[Lemma~3.3]{Liu2024}. This is the mechanism of the counterexample
in \cite[Section~3]{Liu2024}.
\end{remark}

\begin{corollary}\label{cor:PF-no-s}
Suppose $F=|\bar y|\phi(x^0,r,z)$ is a smooth Finsler metric.
Then $F$ is projectively flat in the given coordinates if and
only if, locally,
\[
 \phi(x^0,r,z)=\chi(z)+za(x^0),
\]
where $\chi$ and $a$ are smooth functions of $z$ and $x^0$,
respectively. Equivalently,
$F=|\bar y|\chi(y^0/|\bar y|)+a(x^0)y^0$ is the sum of a Minkowski
norm and the closed one-form $a(x^0)\dd x^0$.
\end{corollary}

\begin{proof}
Since $\phi_s=0$, we have $D\phi=z\phi_{x^0}+\frac sr\phi_r$, hence
$p_3=-\phi_r/r$.
Thus $p_2=p_3=0$ is equivalent to
\[
 \phi_r=0,\qquad z\phi_{x^0z}-\phi_{x^0}=0.
\]
The second equation gives
\[
 \phi_{x^0}=zb(x^0)
\]
for some smooth function $b$. When the domain contains $z=0$,
smoothness ensures that the expressions obtained for $z>0$
and $z<0$ have the same coefficient $b(x^0)$.

Taking a smooth primitive $a$ of $b$ and integrating with
respect to $x^0$, we obtain
\[
 \phi(x^0,r,z)=\chi(z)+za(x^0),
\]
where the integration function $\chi$ is independent of $r$
because $\phi_r=0$.
Conversely, direct substitution gives $p_2=p_3=0$.
\end{proof}

For the general case, put
\begin{equation}\label{eq:characteristic-vars}
 \xi=x^0-zs,\qquad q=r^2-s^2.
\end{equation}
All local reconstruction statements below are made near
$(x^0,r,s,z)=(x^0_*,r_*,0,0)$, $r_*>0$.
The neighborhoods are chosen small enough to contain every
integration segment appearing in the formulas.

\begin{lemma}\label{lem:PF-A}
Let $\phi$ be smooth near $(x^0_*,r_*,0,0)$, where $r_*>0$,
and suppose that $p_2=p_3=0$. Then, locally,
\begin{equation}\label{eq:A-representation}
 \Omega=\Apot(\xi,q,z),\qquad
 \xi=x^0-zs,\quad q=r^2-s^2,
\end{equation}
where $\Apot$ is smooth near $(x^0_*,r_*^2,0)$ and satisfies
\begin{equation}\label{eq:Aconditions}
 2\Apot_{qz}+z\Apot_{\xi\xi}=0,\qquad \Apot_z(\xi,q,0)=0.
\end{equation}
We call $\Apot$ the characteristic potential of $\phi$.
\end{lemma}

\begin{proof}
Differentiating $\Omega=\phi-s\phi_s-z\phi_z$ with respect
to $s$, we obtain
\[
 \Omega_s=-s\phi_{ss}-z\phi_{sz}.
\]
By Proposition~\ref{prop:equivalent-systems},
\[
 \Omega_{x^0}=\phi_{sz},\qquad
 \frac{\Omega_r}{r}=\phi_{ss}.
\]
Consequently, $D\Omega=0$, that is,
\[
 z\Omega_{x^0}+\frac{s}{r}\Omega_r+\Omega_s=0.
\]
The vector field $D$ of \eqref{eq:D} satisfies
\[
 D\xi=Dq=Dz=0,\qquad Ds=1.
\]
Since $r_*>0$, the variables $(\xi,q,s,z)$ form local
coordinates. Thus $D\Omega=0$ gives
$\Omega=\Apot(\xi,q,z)$ for some smooth function $\Apot$.

Taking the $z$-derivative of $\Omega_r/r=\phi_{ss}$
and the $s$-derivative of $\Omega_{x^0}=\phi_{sz}$,
we obtain
\[
 \frac{\Omega_{rz}}r
 =\phi_{ssz}
 =\Omega_{x^0s}.
\]
Using \eqref{eq:A-representation}, the chain rule gives
\[
 \frac{\Omega_{rz}}r
 =2\Apot_{qz}-2s\Apot_{\xi q},\qquad
 \Omega_{x^0s}
 =-z\Apot_{\xi\xi}-2s\Apot_{\xi q}.
\]
Therefore,
\[
 2\Apot_{qz}+z\Apot_{\xi\xi}=0.
\]

Finally,
\[
 \Omega_z=-s\phi_{sz}-z\phi_{zz}
          =-s\Apot_\xi+\Apot_z.
\]
Evaluating at $s=z=0$, where $\xi=x^0$ and $q=r^2$,
we conclude that $\Apot_z(\xi,q,0)=0$.
\end{proof}

\begin{remark}\label{rem:wave}
If $F$ is reversible, then $\phi(x^0,r,-s,-z)=\phi(x^0,r,s,z)$, so
$\Omega$ has the same symmetry. Since $\xi$ and $q$ are invariant
under $(s,z)\mapsto(-s,-z)$, the potential is even in $z$, and we may
write $\Apot(\xi,q,z)=\widetilde\Apot(\xi,q,t)$ with $t=z^2$ and
$\widetilde\Apot$ smooth. Then $\Apot_z(\xi,q,0)=0$ holds
automatically, and the first equation in \eqref{eq:Aconditions}
becomes
\[
 \widetilde\Apot_{\xi\xi}+4\widetilde\Apot_{qt}=0 .
\]
In the variables $a=(q+t)/2$, $b=(q-t)/2$, this is the
$(2+1)$-dimensional wave equation
\[
 \widetilde\Apot_{\xi\xi}+\widetilde\Apot_{aa}=\widetilde\Apot_{bb},
\]
with $b$ playing the role of time.
By Whitney's theorem on even functions (with smooth dependence on the
parameters $\xi,q$), $\widetilde\Apot$ is smooth up to $t=0$ and admits
smooth extensions to $t<0$. However, the equation holds only on
$t\geq0$, that is, on the half-space $a\geq b$. Thus the characteristic
potentials of reversible metrics correspond to solutions of this wave
equation on the region $a\geq b$, and
Theorem~\ref{thm:PF-reconstruction} converts each potential satisfying
\eqref{eq:Aconditions} into a smooth local solution of the reduced
Hamel equations.
\end{remark}

\begin{theorem}\label{thm:PF-reconstruction}
Let $\Apot(\xi,q,z)$ be smooth near $(\xi_*,q_*,0)$,
where $q_*>0$, and satisfy \eqref{eq:Aconditions}. Define
\begin{align}
 B_0(x^0,r,z)
 &=\Apot(x^0,r^2,0)
 -\int_0^z(z-\zeta)
 \frac{\Apot_z(x^0,r^2,\zeta)}{\zeta}\dd\zeta,
 \label{eq:B0}\\
 B_1(x^0,r,z)
 &=\int_0^z\Apot_\xi(x^0,r^2,\zeta)\dd\zeta,
 \label{eq:B1}\\
 \Phi_{\Apot}(x^0,r,s,z)
 &=B_0(x^0,r,z)+sB_1(x^0,r,z)\nonumber\\
 &\quad+2\int_0^s(s-\tau)
 \Apot_q(x^0-z\tau,r^2-\tau^2,z)\dd\tau.
 \label{eq:PhiA}
\end{align}
Then, near $(\xi_*,\sqrt{q_*},0,0)$, the function
$\Phi_{\Apot}$ is smooth, satisfies $p_2=p_3=0$, and has
characteristic potential $\Apot$.
Moreover, every smooth solution with this potential is locally of the form
\begin{equation}\label{eq:PF-family}
 \phi=\Phi_{\Apot}
 +z\psi_{x^0}(x^0,r)+\frac{s}{r}\psi_r(x^0,r),
\end{equation}
where $\psi$ is smooth. Conversely, every function of this form
satisfies $p_2=p_3=0$ and has characteristic potential $\Apot$.
As $\Apot$ ranges over the smooth solutions of
\eqref{eq:Aconditions}, this construction gives every smooth
local solution near a point with $r>0$ and $s=z=0$.
\end{theorem}

\begin{proof}
Since $\Apot_z(\xi,q,0)=0$, we have
\[
 \frac{\Apot_z(\xi,q,z)}{z}
 =\int_0^1\Apot_{zz}(\xi,q,tz)\dd t.
\]
Thus the quotient extends smoothly to $z=0$, and the functions
in \eqref{eq:B0}--\eqref{eq:PhiA} are smooth.

Put $\Phi=\Phi_{\Apot}$ and
\[
 \widehat\Apot(x^0,r,s,z)
 =\Apot(x^0-zs,r^2-s^2,z).
\]
Differentiating \eqref{eq:PhiA}, we obtain
\[
 \Phi_{ss}
 =2\Apot_q(x^0-zs,r^2-s^2,z)
 =\frac{\widehat\Apot_r}{r}.
\]
The equation $2\Apot_{qz}+z\Apot_{\xi\xi}=0$ gives
\[
 \frac{\widehat\Apot_{rz}}{r}
 =\widehat\Apot_{x^0s}.
\]
Consequently,
\[
 \partial_s\bigl(\Phi_{sz}-\widehat\Apot_{x^0}\bigr)=0.
\]
At $s=0$, \eqref{eq:B1} gives
\[
 \Phi_{sz}=(B_1)_z
 =\Apot_\xi(x^0,r^2,z)
 =\widehat\Apot_{x^0}.
\]
Hence $\Phi_{sz}=\widehat\Apot_{x^0}$.

Since $D\widehat\Apot=0$, it follows that
\begin{align*}
 \partial_s\bigl(\Phi-s\Phi_s-z\Phi_z-\widehat\Apot\bigr)
 &=-s\Phi_{ss}-z\Phi_{sz}-\widehat\Apot_s\\
 &=-\frac{s}{r}\widehat\Apot_r
   -z\widehat\Apot_{x^0}-\widehat\Apot_s\\
 &=0.
\end{align*}
Direct differentiation of \eqref{eq:B0} gives
\[
 B_0-z(B_0)_z=\Apot(x^0,r^2,z).
\]
Thus the preceding expression vanishes at $s=0$, proving
\[
 \Phi-s\Phi_s-z\Phi_z=\widehat\Apot.
\]

To verify the remaining compatibility in \eqref{eq:PF-omega},
put
\[
 C=\frac{\Phi_{rz}}{r}-\Phi_{x^0s}.
\]
Using the identities above, we obtain
\[
 C_s
 =\frac{\widehat\Apot_{x^0r}}{r}
  -\frac{\widehat\Apot_{rx^0}}{r}
 =0.
\]
At $s=0$, \eqref{eq:B0} and \eqref{eq:B1} give
\[
 C=-2\int_0^z
 \frac{\Apot_{qz}(x^0,r^2,\zeta)}{\zeta}\dd\zeta
 -\int_0^z\Apot_{\xi\xi}(x^0,r^2,\zeta)\dd\zeta
 =0,
\]
by \eqref{eq:Aconditions}, with the quotient understood by
smooth extension at $\zeta=0$.
Therefore $C=0$, and
Proposition~\ref{prop:equivalent-systems} yields $p_2=p_3=0$.

Multiplying the additional term in \eqref{eq:PF-family} by
$u=|\bar y|$, we obtain $d\psi_x(y)$, where
\[
 d\psi
 =\psi_{x^0}\,dx^0
 +\frac{\psi_r}{r}\sum_i x^i\,dx^i.
\]
The Hamel residual of a one-form $\theta=\theta_Ay^A$ is
\[
 H_B(\theta)
 =(\partial_A\theta_B-\partial_B\theta_A)y^A.
\]
It vanishes for $d\psi$, so the additional term preserves
$p_2=p_3=0$. Since this term is linear in $(s,z)$, it also
preserves $\Omega$ (compare \cite[Lemmas~3.2, 3.3 and~4.1]{Liu2024}).

Conversely, let $\phi$ be a solution with characteristic
potential $\Apot$, and put $\delta=\phi-\Phi_{\Apot}$.
By \eqref{eq:PF-omega},
\[
 \delta_{ss}=\delta_{sz}=0,\qquad
 \delta-s\delta_s-z\delta_z=0.
\]
The first two identities give
\[
 \delta=b(x^0,r)s+h(x^0,r,z).
\]
The last identity becomes $h-zh_z=0$. Smoothness at $z=0$
then gives $h=a(x^0,r)z$. Hence
\[
 \delta=a(x^0,r)z+b(x^0,r)s.
\]
The remaining compatibility in \eqref{eq:PF-omega} yields
\[
 a_r=rb_{x^0}.
\]
Thus the one-form $a\,dx^0+rb\,dr$ is closed and, on a smaller
rectangular base domain, exact. Therefore
\[
 a=\psi_{x^0},\qquad b=\frac{\psi_r}{r}
\]
for some smooth function $\psi$, proving \eqref{eq:PF-family}.
Finally, Lemma~\ref{lem:PF-A} shows that every smooth local
solution near such a point determines a potential satisfying
\eqref{eq:Aconditions}.
\end{proof}

\begin{remark}\label{rem:conic}
The reconstruction gives smooth local solutions of
$p_2=p_3=0$ near a point with $r>0$ and $s=z=0$.
It imposes neither positivity nor strong convexity and
does not guarantee extension across $r=0$, $|s|=r$,
or $\bar y=0$.

On any open conic domain where $F=|\bar y|\phi$ is positive
and strongly convex, it defines a conic Finsler metric
in the sense of \cite{JavaloyesSanchez2014}.
Hamel's equations then imply that $F$ is projectively flat.
\end{remark}

\begin{proposition}\label{prop:potential-positive}
For a solution of \eqref{eq:corrected} with
$\Omega=\Apot(\xi,q,z)$, one has
\[
 \phi_{ss}=2\Apot_q,\qquad \phi_{sz}=\Apot_\xi,\qquad
 \phi_{zz}=-\frac{\Apot_z}{z},
\]
where $\Apot_z$ denotes differentiation with $(\xi,q)$ fixed and
$-\Apot_z/z$ is extended at $z=0$ by
$-\Apot_{zz}(\xi,q,0)$. Moreover,
\begin{equation}\label{eq:Lambda-A}
 \Lambda=-\frac{\Apot_z}{z}(\Apot+2q\Apot_q)-q\Apot_\xi^2.
\end{equation}
In particular, if $\phi>0$, then $F$ is strongly convex on the regular
reduced region if and only if $-\Apot_z/z>0$, $\Lambda>0$ and,
when $n\geq3$, $\Apot>0$.
\end{proposition}
\begin{proof}
The first two identities follow from \eqref{eq:PF-omega}.
On the one hand, $\Omega_z=-s\phi_{sz}-z\phi_{zz}$.
On the other hand, differentiating the composition gives
$\Omega_z=-s\Apot_\xi+\Apot_z$. Since $\phi_{sz}=\Apot_\xi$,
we obtain $\phi_{zz}=-\Apot_z/z$, including $z=0$ by smooth
extension. Substitution gives \eqref{eq:Lambda-A},
and Lemma~\ref{lem:convexity} gives the convexity criterion.
\end{proof}

\begin{example}\label{ex:local}
Let $\delta>0$ and
\[
 \Apot(\xi,q,z)=1+\delta\bigl(2\xi^2-qz^2\bigr).
\]
Then $2\Apot_{qz}+z\Apot_{\xi\xi}=-4\delta z+4\delta z=0$ and
$\Apot_z=-2\delta qz$ vanishes at $z=0$. Formulas
\eqref{eq:B0}--\eqref{eq:PhiA} give
$B_0=1+2\delta(x^0)^2+\delta r^2z^2$, $B_1=4\delta x^0z$ and
\[
 \Phi_{\Apot}=1+\delta\bigl(2(x^0)^2+r^2z^2+4x^0sz-s^2z^2\bigr).
\]
A direct check confirms $\Omega=\Apot$ and \eqref{eq:PF-omega}.
In the original variables,
\[
 F=|\bar y|+\delta\left(2(x^0)^2|\bar y|
   +\frac{r^2(y^0)^2}{|\bar y|}
   +\frac{4x^0y^0\ip{\bar x}{\bar y}}{|\bar y|}
   -\frac{(y^0)^2\ip{\bar x}{\bar y}^2}{|\bar y|^3}\right).
\]
By Proposition~\ref{prop:potential-positive},
$\phi_{zz}=2\delta q$ and
\[
 \Lambda=2\delta q\bigl(\Apot-2\delta qz^2\bigr)-16\delta^2q\xi^2 ,
\]
so that $\Lambda=2\delta r^2\bigl(1-6\delta(x^0)^2\bigr)$ at $s=z=0$.
Hence, near every reduced point $(x^0,r,0,0)$ with $r>0$ and
$6\delta(x^0)^2<1$, $F$ is a reversible, cylindrically symmetric,
projectively flat conic Finsler metric \cite{JavaloyesSanchez2014} on
an open $\OO(n)$-invariant conic domain of tangent vectors.
It is only local: it is not defined at $\bar y=0$.
Since $\Apot(\xi,q,0)=1+2\delta\xi^2$ does not depend on $q$, this
potential is not of the form \eqref{eq:profile-A} below for any $c$
and $f$. In the reversible normal form of Remark~\ref{rem:wave},
$\widetilde\Apot=1+\delta(2\xi^2-qt)$ is a polynomial solution of
$\widetilde\Apot_{\xi\xi}+4\widetilde\Apot_{qt}=0$.
\end{example}

For the following metrics, smoothness and strong convexity are
verified on $TM_0$ in the original variables.

\section{Characteristic potentials of a spherical family}
\label{sec:construction}

We now apply the preceding formulas to the projectively flat spherically symmetric metrics of \cite{Zhou2012,HuangMo2013}
and to their pullbacks under the linear maps
\[
 L(x^0,\bar x)=(\lambda_0x^0,\lambda_1\bar x),\qquad
 \lambda_0,\lambda_1>0,
\]
which rescale the axial and transverse directions differently when
$\lambda_0\neq\lambda_1$. A direct computation gives
$L^*F_{1,f}=\lambda_1F_{c,\tilde f}$ with $c=\lambda_0/\lambda_1$ and
$\tilde f(w)=f(\lambda_1^2w)$, in the notation of
Theorem~\ref{thm:construction} below. Hence it suffices to consider
the one-parameter family $T_c(x^0,\bar x)=(cx^0,\bar x)$, $c>0$.
For $c=1$ the metrics in
Theorem~\ref{thm:construction} are spherically symmetric and known;
see Remark~\ref{rem:origin}. The new information is the
characteristic potential \eqref{eq:profile-A} and the convexity
criterion, which only requires $f>0$ and $f+2wf'>0$.

\begin{theorem}\label{thm:construction}
Let $c>0$, $0<R\leq\infty$, and let $f$ be smooth on an
open neighborhood of $[0,R)$. Suppose that
\begin{equation}\label{eq:f-hyp}
 f(w)>0\qquad\text{and}\qquad f(w)+2wf'(w)>0
 \qquad\text{for every }w\in[0,R).
\end{equation}
Define
\begin{equation}\label{eq:geometric-vars}
 \alpha_c=\sqrt{c^2(y^0)^2+|\bar y|^2},\qquad
 X=c^2(x^0)^2+r^2,\qquad
 \sigma=\frac{c^2x^0y^0+\ip{\bar x}{\bar y}}{\alpha_c}.
\end{equation}
For $X<R$, put
\begin{equation}\label{eq:H}
 H(X,\sigma)=f(X-\sigma^2)
 +2\sigma\int_0^\sigma f'(X-\tau^2)\dd\tau.
\end{equation}
Then
\begin{equation}\label{eq:family}
 F_{c,f}(x,y)=\alpha_c H(X,\sigma)
\end{equation}
is a smooth reversible cylindrically symmetric projectively flat
Finsler metric on the entire slit tangent bundle. Its characteristic
potential is
\begin{equation}\label{eq:profile-A}
 \Apot(\xi,q,z)=\frac{f(w)}\varrho,\qquad
 \varrho=\sqrt{1+c^2z^2},\qquad
 w=q+\frac{c^2\xi^2}{\varrho^2}.
\end{equation}
Moreover,
\begin{equation}\label{eq:profile-positive}
 \Omega=\frac{f(w)}\varrho>0,\qquad
 \Lambda=\frac{c^2}{\varrho^4}f(w)\bigl(f(w)+2wf'(w)\bigr)>0.
\end{equation}
In particular, the condition $f+2wf'>0$ in \eqref{eq:f-hyp} is
equivalent to $\Lambda>0$.
\end{theorem}
\begin{proof}
First take $c=1$, so that $\alpha_1=|y|$, $X=|x|^2$ and
$\sigma=\ip{x}{y}/|y|$. Cauchy--Schwarz gives
$0\leq X-\sigma^2\leq X<R$. Formula \eqref{eq:H} can be written as
\[
 H(X,\sigma)=f(X-\sigma^2)
 +2\sigma^2\int_0^1 f'(X-\sigma^2t^2)\dd t.
\]
Thus $H$ is smooth and even in $\sigma$. Since
$\alpha_1=|y|>0$ for $y\ne0$, $F$ is smooth on the entire
slit tangent bundle, reversible, and positively homogeneous.

Differentiation gives
\begin{equation}\label{eq:Hderivatives}
 H_\sigma=2\int_0^\sigma f'(X-\tau^2)\dd\tau,\quad
 H_{\sigma\sigma}=2f'(X-\sigma^2),\quad
 H-\sigma H_\sigma=f(X-\sigma^2).
\end{equation}
Fix $x$ and $y\ne0$, and put $w=X-\sigma^2$. For $V\perp y$,
the fiber Hessian is
\[
 D_y^2F(V,V)=\frac1{|y|}
 \left[f(w)|V|^2+2f'(w)\ip{x}{V}^2\right].
\]
Since $\ip{x}{V}=\ip{x-\sigma y/|y|}{V}$ and
$|x-\sigma y/|y||^2=w$, we have $0\leq\ip{x}{V}^2\leq w|V|^2$.
Hence
\[
 D_y^2F(V,V)\geq\frac1{|y|}\min\bigl\{f(w),\,f(w)+2wf'(w)\bigr\}|V|^2>0
 \qquad(V\ne0),
\]
by \eqref{eq:f-hyp}. Homogeneity gives $D_y^2F(y,\cdot)=0$, so
$D_y^2F$ is positive semidefinite with kernel $\R y$.

Next we show that $F>0$. If $x=0$, then $F(0,y)=f(0)|y|>0$. Let
$x\neq0$, choose a unit vector $e\perp x$, and put
$y(\theta)=\cos\theta\,x/|x|+\sin\theta\,e$ and
$k(\theta)=F(x,y(\theta))=H(X,|x|\cos\theta)$. Since $F$ is
homogeneous of degree one and $y''=-y$, we have
$k+k''=D_y^2F(y',y')>0$ with $y'\perp y$. The function
$m(\theta)=k(\theta)-k(\pi/2)\sin\theta$ satisfies $m+m''>0$,
$m(\pi/2)=0$ and $m'(\pi/2)=-|x|H_\sigma(X,0)=0$. Hence
\[
 m(\theta)=\int_{\pi/2}^{\theta}\sin(\theta-t)\,(m+m'')(t)\dd t>0
 \qquad(\pi/2<\theta\leq\pi),
\]
and therefore $k(\theta)>k(\pi/2)\sin\theta\geq0$ on $[\pi/2,\pi]$,
where $k(\pi/2)=f(X)>0$. Thus $H(X,\sigma)>0$ for
$-|x|\leq\sigma\leq0$, and by evenness for all $|\sigma|\leq|x|$.
Every pair $(x,y)$ with $y\neq0$ is a positive multiple of one of these
$(x,y(\theta))$ for a suitable choice of $e$, so $F>0$ on $TM_0$.
Consequently $g_y=d_yF\otimes d_yF+F D_y^2F$ is positive
definite for every nonzero $y$, including $\bar y=0$.

Taking the $X$-derivative of the last identity in
\eqref{eq:Hderivatives}, we obtain
\[
 H_X-\sigma H_{X\sigma}=f'(X-\sigma^2)
 =\frac12H_{\sigma\sigma}.
\]
For $F=|y|H(X,\sigma)$, direct differentiation gives
\[
 F_{x^Ay^B}y^A-F_{x^B}
 =\bigl(H_{\sigma\sigma}+2\sigma H_{X\sigma}-2H_X\bigr)
   \bigl(|y|x_B-\sigma y_B\bigr)=0.
\]
Thus $F$ is projectively flat. For general $c>0$, formula
\eqref{eq:family} is its pullback by the invertible linear map
$(x^0,\bar x)\mapsto(cx^0,\bar x)$, which preserves straight
lines and commutes with the spatial orthogonal action.

In reduced variables,
\[
 \phi=\varrho H(X,\sigma),\qquad
 \sigma=\frac{c^2x^0z+s}{\varrho},\qquad
 X-\sigma^2=q+\frac{c^2\xi^2}{\varrho^2}=w.
\]
The Euler expression gives
$\Omega=(H-\sigma H_\sigma)/\varrho=f(w)/\varrho$.
For \eqref{eq:profile-A}, one has
\[
 \Apot_q=\frac{f'}\varrho,\qquad
 \Apot_\xi=\frac{2c^2\xi f'}{\varrho^3},\qquad
 -\frac{\Apot_z}{z}=\frac{c^2f}{\varrho^3}
 +\frac{2c^4\xi^2f'}{\varrho^5}.
\]
These formulas verify $2\Apot_{qz}+z\Apot_{\xi\xi}=0$ and, after
substitution in \eqref{eq:Lambda-A}, give
\eqref{eq:profile-positive}.
\end{proof}

\begin{remark}\label{rem:origin}
For $c=1$, $F_{1,f}$ is spherically symmetric and belongs to the
projectively flat family classified in \cite{Zhou2012,HuangMo2013}.
Formula \eqref{eq:H} also follows from
\cite[Remark~1, equations~(3.22)--(3.24)]{Solorzano2023}
by choosing the profile $2f'$, the constant $f(0)$, and the
odd term equal to zero; it is the smooth even solution of
$H-\sigma H_\sigma=f(X-\sigma^2)$.
For $c>0$, the map $T_c(x^0,\bar x)=(cx^0,\bar x)$ gives
$F_{c,f}=T_c^*F_{1,f}$. Thus $F_{c,f}$ is isometric to a spherically
symmetric metric, and for each single profile the parameter $c$ does
not produce new metrics up to linear changes of coordinates.
It does change the characteristic potential \eqref{eq:profile-A},
and this is what Example~\ref{ex:sum} exploits.
\end{remark}

\begin{proposition}\label{prop:nonriemannian}
At a point $x$ with $X>0$, the metric in
Theorem~\ref{thm:construction} is non-Riemannian if $f^{-2}$
is not affine on $[0,X]$. In particular, this holds if
$ff''-3(f')^2$ is not identically zero on that interval.
In that case the norm at that point is not Randers.
\end{proposition}
\begin{proof}
Work in the linearly transformed coordinates with $c=1$, in which
$F$ is invariant under the full group $\OO(n+1)$ acting on $(x,y)$.
Fix $x\neq0$. The stabilizer of $x$ acts on $x^\perp$, whose dimension
$n\geq2$, as the full orthogonal group. Therefore a Riemannian norm at
$x$ invariant under this stabilizer has the form
$a|y|^2+b\ip{x}{y}^2$, that is,
$H(X,\sigma)^2=C+D\sigma^2$, where $C>0$ and $D$ depend on the
fixed point. Therefore
\[
 f(X-\sigma^2)=H-\sigma H_\sigma
 =\frac{C}{\sqrt{C+D\sigma^2}}.
\]
Putting $w=X-\sigma^2$ gives
$f(w)^{-2}=[C+D(X-w)]/C^2$, which is affine.
Also,
\[
 (f^{-2})''=-2f^{-4}\bigl(ff''-3(f')^2\bigr).
\]
Finally, a reversible Randers metric $\alpha+\beta$ has
$\beta=0$ and is Riemannian. The claim follows. At the origin, the
norm is $F_{c,f}(0,y)=f(0)\alpha_c(y)$, which is Riemannian.
\end{proof}

\begin{remark}\label{rem:sphere}
The conditions in \eqref{eq:f-hyp} also allow decreasing
functions $f$. For instance,
$f(w)=(1+\eta w)^{-p}$ with $\eta>0$ satisfies
\[
 f+2wf'=(1+\eta w)^{-p-1}\bigl(1+(1-2p)\eta w\bigr)>0
\]
for all $w\geq0$ if and only if $0<p\leq1/2$. For $p=1/2$,
$f^{-2}=1+\eta w$ is affine, and the proof of
Proposition~\ref{prop:nonriemannian} gives, for $c=1$,
\[
 F_{1,f}^2=\frac{|y|^2(1+\eta|x|^2)-\eta\ip{x}{y}^2}{(1+\eta|x|^2)^2},
\]
the round metric of constant sectional curvature $\eta$ in
gnomonic coordinates. For $0<p<1/2$, $f^{-2}$ is not affine, and
these metrics are non-Riemannian at every $x$ with $X>0$.
The previous hypothesis $f'\geq0$ excluded all these profiles.
\end{remark}

\section{Examples}\label{sec:examples}
We use the notation in \eqref{eq:geometric-vars}. Each metric
is given in the original variables, so its regularity includes
all nonzero tangent directions. By Remark~\ref{rem:origin},
the metrics in Examples~\ref{ex:linear}--\ref{ex:exponential} are
linear pullbacks of known spherically symmetric metrics; we record
them for their explicit characteristic potentials.
Example~\ref{ex:sum} is genuinely cylindrical.

\begin{example}\label{ex:linear}
Taking $f(w)=1+\eta w$, where $c,\eta>0$, formula \eqref{eq:H} gives
\[
 H(X,\sigma)=1+\eta(X+\sigma^2).
\]
Thus the metric
\begin{equation}\label{eq:linear-metric}
 F=(1+\eta X)\alpha_c
 +\eta\frac{(c^2x^0y^0+\ip{\bar x}{\bar y})^2}{\alpha_c}
\end{equation}
is defined on $\R^{n+1}$, with
\[
 \Apot=\frac{1+\eta w}{\varrho},\qquad
 \Lambda=\frac{c^2}{\varrho^4}(1+\eta w)(1+3\eta w)>0.
\]
Theorem~\ref{thm:construction} proves positivity, smoothness,
and strong convexity on the full slit tangent bundle.
Since $ff''-3(f')^2=-3\eta^2$, the norm is non-Riemannian
at every point with $X>0$, and is not Randers there.
\end{example}

\begin{example}\label{ex:quadratic}
Taking $f(w)=1+\eta w^2$, where $c,\eta>0$, we obtain
\[
 H(X,\sigma)=1+\eta X^2+2\eta X\sigma^2-\frac{\eta}{3}\sigma^4.
\]
Consequently,
\begin{equation}\label{eq:quadratic-metric}
 F=\alpha_c\left(1+\eta X^2+2\eta X\sigma^2
                      -\frac{\eta}{3}\sigma^4\right)
\end{equation}
is a projectively flat Finsler metric on $\R^{n+1}$. Indeed,
$f>0$ and $f+2wf'=1+5\eta w^2>0$, and
\[
 \Apot=\frac{1+\eta w^2}{\varrho},\qquad
 \Lambda=\frac{c^2}{\varrho^4}(1+\eta w^2)(1+5\eta w^2)>0.
\]
The negative coefficient of $\sigma^4$ does not affect the
positivity or strong convexity, by Theorem~\ref{thm:construction}.
The function $ff''-3(f')^2=2\eta-10\eta^2w^2$ is not
identically zero on any interval of positive length.
Thus the norm is non-Riemannian and non-Randers for $X>0$.
\end{example}

\begin{example}\label{ex:exponential}
Taking $f(w)=e^{\eta w}$, where $c,\eta>0$, formula \eqref{eq:H} gives
\[
 H(X,\sigma)=e^{\eta(X-\sigma^2)}
 +\sqrt{\pi\eta}\,\sigma e^{\eta X}\erf(\sqrt\eta\,\sigma),
\]
where $\erf(u)=2\pi^{-1/2}\int_0^u e^{-t^2}\dd t$.
For $c=1$, this is the metric in
\cite[Theorem~1.2]{HuangMo2013} with $g=0$ and $\lambda=\eta$.
Hence
\begin{equation}\label{eq:exponential-metric}
 F=\alpha_c\left[e^{\eta(X-\sigma^2)}
 +\sqrt{\pi\eta}\,\sigma e^{\eta X}\erf(\sqrt\eta\,\sigma)\right]
\end{equation}
is a smooth strongly convex projectively flat Finsler metric
on $\R^{n+1}$. Its reduced expressions are
\[
 \Apot=\frac{e^{\eta w}}\varrho,\qquad
 \Lambda=\frac{c^2}{\varrho^4}e^{2\eta w}(1+2\eta w)>0.
\]
Since $ff''-3(f')^2=-2\eta^2e^{2\eta w}$, it is non-Riemannian
and non-Randers at every point with $X>0$.
\end{example}

\begin{proposition}\label{prop:sums}
A positive finite sum of smooth strongly convex projectively flat
Finsler metrics in the same coordinates is again a smooth strongly
convex projectively flat Finsler metric. Cylindrical symmetry and
reversibility are preserved when they hold for every summand.
The characteristic potential of the sum is the corresponding sum of
characteristic potentials.
\end{proposition}
\begin{proof}
Let $F=\sum_j\lambda_jF_j$, $\lambda_j>0$.
The fiber Hessians $D_y^2F_j$ are nonnegative with the same
radial kernel. Their sum has the same property. Positivity of $F$
therefore gives positive definiteness of
$d_yF\otimes d_yF+F D_y^2F$. Hamel's equations and the map
$\phi\mapsto\Omega$ are linear.
The remaining properties follow from the definition.
\end{proof}

\begin{example}\label{ex:sum}
Let $\eta,\mu>0$. Consider the sum
\[
 F=F_{1,f_1}+F_{2,f_2},\qquad
 f_1(w)=1+\eta w,\quad f_2(w)=1+\mu w^2.
\]
Proposition~\ref{prop:sums} gives a smooth reversible cylindrically
symmetric projectively flat Finsler metric on $\R^{n+1}$.
Its potential is
\[
 \Apot=\frac{1+\eta w_1}{\sqrt{1+z^2}}
   +\frac{1+\mu w_2^2}{\sqrt{1+4z^2}},\qquad
 w_1=q+\frac{\xi^2}{1+z^2},\quad
 w_2=q+\frac{4\xi^2}{1+4z^2}.
\]
At the origin,
\[
 F(0,y)=\sqrt{(y^0)^2+|\bar y|^2}
       +\sqrt{4(y^0)^2+|\bar y|^2}.
\]
This is non-Riemannian. Otherwise, squaring would imply that the
product of two nonproportional positive definite quadratic forms
is the square of a quadratic polynomial. The two forms are
irreducible over $\R$, so unique factorization excludes that
possibility. Reversibility also excludes a Randers representation.

Consequently, $F$ is not the pullback of a spherically symmetric
metric by a linear isomorphism $L$ of $\R^{n+1}$. Indeed, if
$G(x,y)=|y|\varphi\bigl(|x|,\ip{x}{y}/|y|\bigr)$, then
$(L^*G)(0,y)=G(0,Ly)=\varphi(0,0)|Ly|$ is Riemannian.
In contrast with Examples~\ref{ex:linear}--\ref{ex:exponential},
this metric is genuinely cylindrically symmetric.
\end{example}

\begin{remark}
Other choices in Theorem~\ref{thm:construction} are
\[
 1+\eta w^3,\qquad (1+\eta w)^p,\qquad
 1+\log(1+\eta w),\qquad 1+\frac{\eta w}{1+w},\qquad
 (1+\eta w)^{-p'} ,
\]
where $\eta,p>0$ and $0<p'\leq1/2$ (see Remark~\ref{rem:sphere}).
The choice $(1-\eta w)^{-p}$ gives metrics on $X<1/\eta$; its case
$p=1/2$ is the Klein model of hyperbolic space.
The corresponding characteristic potentials are obtained by
substituting these profiles in \eqref{eq:profile-A}, and positive sums
of the corresponding $F_{c_j,f_j}$ with distinct $c_j$ give further
cylindrically symmetric examples as in Example~\ref{ex:sum}.
\end{remark}

\begin{remark}
For comparison with the explicit separated family in
\cite[Theorem~3]{Solorzano2023}, every member of that family
has $\phi_{rz}=\phi_{x^0s}=0$. For the present construction,
\[
 \phi_{sz}=\frac{2c^2\xi f'(w)}{\varrho^3},\qquad
 \phi_{x^0sz}\big|_{x^0=s=z=0}=2c^2f'(r^2).
\]
Hence, when $f'(r^2)\ne0$, these reduced functions are outside
that separated formula near the indicated point. The comparison is
made in the coordinates of \eqref{def:F}.
We note that the separated family of \cite[Theorem~3]{Solorzano2023}
remains projectively flat under the characterization of
\cite{Liu2024}: it satisfies the first two equations of
\eqref{eq:PF-omega} by construction, and the third one holds
trivially because $\phi_{rz}=\phi_{x^0s}=0$
(Remark~\ref{rem:missing}).
\end{remark}

\end{document}